\documentclass[11pt]{article}
\usepackage{amsfonts,epsf,amsmath,amssymb}
\usepackage{tikz}

\newtheorem{theorem}{\bf Theorem}[section]
\newtheorem{corollary}[theorem]{\bf Corollary}

\newtheorem{proposition}[theorem]{\bf Proposition}

\newcommand{\vertex}{\node[vertex]}
\tikzstyle{vertex}=[circle, draw, inner sep=0pt, minimum size=6pt]
\newcommand{\pch}{\chi_{\rho}}
\newcommand{\proof}{\noindent{\bf Proof.\ }}
\newcommand{\qed}{\hfill $\square$ \bigskip}

\begin{document}

\title{\bf Packing chromatic number under local changes in a graph}

\author{
 \phantom{xxxxx} \and Bo\v{s}tjan Bre\v{s}ar $^{a,b}$  \and Sandi Klav\v zar $^{a,b,c}$ \and \phantom{xxx}
\and Douglas F. Rall $^{d}$ \and Kirsti Wash $^e$}

\date{}

\maketitle

\begin{center}
$^a$ Faculty of Natural Sciences and Mathematics, University of Maribor, Slovenia

$^b$ Institute of Mathematics, Physics and Mechanics, Ljubljana, Slovenia

$^c$ Faculty of Mathematics and Physics, University of Ljubljana, Slovenia

$^d$ Department of Mathematics, Furman University, Greenville, SC, USA

$^e$ Department of Mathematics, Trinity College, Hartford, CT, USA
\end{center}

\begin{abstract}
The packing chromatic number $\chi_{\rho}(G)$ of a graph $G$ is the
smallest integer $k$ such that there exists a $k$-vertex coloring of
$G$ in which any two vertices receiving color $i$ are at distance at
least $i+1$. It is proved that in the class of subcubic graphs the
packing chromatic number is bigger than $13$, thus answering an open
problem from [Gastineau, Togni, $S$-packing colorings of cubic
graphs, Discrete Math.\ 339 (2016) 2461--2470]. In addition, the
packing chromatic number is investigated with respect to several
local operations. In particular, if $S_e(G)$ is the
graph obtained from a graph $G$ by subdividing its edge $e$, then
$\left\lfloor \pch(G)/2 \right\rfloor +1 \le \chi_{\rho}(S_e(G)) \le
\pch(G)+1$.
\end{abstract}


\medskip\noindent
{\bf AMS Subj.\ Class:} 05C70, 05C15, 05C12

\baselineskip16pt

\section{Introduction}
Many variations of the classical graph coloring have been
introduced, several of which involve graph distance, which as a
condition is usually imposed on the vertices that are given the same
color. In this paper we study packing colorings defined as follows.
The {\em packing chromatic number} $\chi_{\rho}(G)$ of $G$ is the
smallest integer $k$ such that $V(G)$ can be partitioned into
subsets $X_1, \ldots, X_k$, where $X_i$ induces an $i$-packing;
that is, vertices of $X_i$ are pairwise at distance more than $i$.
Equivalently, a {\em $k$-packing coloring} of $G$ is a function
$c:V(G)\rightarrow [k]$, where $[k]=\{1,\ldots ,k\}$, such that if
$c(u) = c(v) = i$, then $d_G(u,v) > i$, where $d_G(u,v)$ is the
usual shortest-path distance between $u$ and $v$ in $G$. We
mention that in distance-$k$ colorings $V(G)$ is partitioned into
$k$-packings.

The concept of the packing chromatic number was introduced
in~\cite{goddard-2008} and given the name in~\cite{bresar-2007}. The
problem intuitively appears more difficult
than the standard coloring problem. Indeed, the packing chromatic
number is intrinsically more difficult due to the fact that determining
$\chi_\rho$ is NP-complete even when restricted to
trees~\cite{fiala-2010}. On the other hand, Argiroffo et al.\
discovered that the packing coloring problem is solvable in
polynomial time for several nontrivial classes of
graphs~\cite{argiroffo-2014}. In addition, the packing chromatic number was
studied on hypercubes~\cite{goddard-2008,torres-2015},
Cartesian product graphs~\cite{jacobs-2013,shao-2015}, and distance
graphs~\cite{ekstein-2014,togni-2014}.

In the seminal paper~\cite{goddard-2008} the following problem was
posed: does there exist an absolute constant $M$, such that
$\pch(G)\le M$ holds for any subcubic graph $G$. (Recall that a graph is {\em subcubic}, if
its largest degree is bounded by $3$.) This problem led to
a lot of research but remains unsolved at the present. In
particular, the packing chromatic number of the infinite hexagonal
lattice is 7 (the upper bound being established
in~\cite{fiala-2009}, the lower bound in~\cite{korze-2014}), hence
the packing chromatic number of any subgraph of the hexagonal
lattice is bounded by $7$. The same bound also holds for subcubic
trees as follows from a result of Sloper~\cite{sloper-2004}. For the
(subcubic) family of base-3 Sierpi\' nski graphs the packing
chromatic number was bounded by 9 in~\cite{bkr-2016+}. The exact 
value of the packing chromatic of some additional
subcubic graphs  was determined in~\cite{ds2016}. Very
recently, Gastineau and Togni~\cite{gt-2016} found a cubic graph
with  packing chromatic number equal to 13 and posed an open
problem which intrigued us: does there exist a cubic graph with
packing chromatic number larger than 13?

We proceed as follows.  In the next section we prove that the answer
to the above question is positive. More precisely, we construct a cubic
graph on 78 vertices with  packing chromatic number at least
$14$. A key technique in the related proof is  edge subdivision.
We hence give a closer look at this operation with respect to its
effect on the packing chromatic number. In particular, the packing
chromatic number does not increase by more than 1 when an edge of a
graph is subdivided, but can decrease by at least 2. In addition, we
prove that the lower bound for the packing chromatic number of an
edge-subdivided graph is bigger than half of the packing chromatic
number of the original graph. Then, in Section~\ref{sec:deletion},
we investigate the effect on the packing chromatic number of the
following local operations: a vertex deletion, an edge deletion, and
an edge contraction. In particular, we demonstrate that the
difference $\pch(G)-\pch(G-e)$ can be arbitrarily large.

\section{Edge subdivision}
\label{sec:subdivision}

In this section we consider the packing chromatic number with
respect to the edge-subdivision operation. If $e$ is an edge of a
graph $G$, then let $S_e(G)$ denote the graph obtained from $G$ by
subdividing the edge $e$. The graph obtained from $G$ by subdividing
all its edges is denoted $S(G)$.

The following theorem is the key for the answer of the above
mentioned question of Gastineau and Togni.

\begin{theorem}
\label{thm:diameter} Suppose that there exists a constant $M$ such
that $\pch(H)\le M$ holds for any subcubic graph $H$. If $G$ is a
subcubic graph such that $\pch(G)=M$, then either $\pch(S_e(G))\le
M-2$ for any $e\in E(G)$, or ${\rm diam}(G)\ge \lceil \frac {M}{2}
\rceil -2$.
\end{theorem}

\proof Let $G$ be a subcubic graph such that $\pch(G)=M$, where $\pch(H)\le M$
for every subcubic graph $H$. If $\pch(S_e(G))\le M-2$ holds for any $e\in E(G)$,
there is nothing to be proved. Hence assume that there exists an
edge $e\in E(G)$ such that $\pch(S_e(G))\ge M-1$. Let $G'$ be the
graph obtained from $G$ by subdividing the edge $e$, and let $x'$ be
the new vertex. Let $G''$ be a copy of $G'$, with $x''$ playing the
role of $x'$. Let now $\widehat{G}$ be the graph obtained from the
disjoint union of $G'$ and $G''$ by connecting $x'$ with $x''$.

Note first that $\widehat{G}$ is a subcubic graph, and hence by the
theorem's assumption, $\pch(\widehat{G})\le M$. Let $c$ be an
arbitrary optimal packing coloring of $\widehat{G}$. Because
$c$ restricted to $G'$ (resp.\ $G''$) is a packing coloring of
$G'=S_e(G)$ (resp.\ $G''$), $c$ uses at least $M-1$ colors.
We claim that ${\rm diam}(\widehat{G})\ge M-1$. If $c$ colors a
vertex $u'$ of $G'$ and a vertex $u''$ of $G''$ by the color $M$,
then $d_{\widehat{G}}(u',u'')> M$, and the claim follows. Otherwise,
we assume that $c$ restricted to $G'$ does not use the color $M$. If also
$G''$ does not use color $M$, then since $\pch(G')\ge M-1$ and
$\pch(G'')\ge M-1$, there exist vertices $v',v''$ in $G'$, resp.\
$G''$, with $c(v')=c(v'')=M-1$, and consequently, ${\rm
diam}(\widehat{G})>M-1$ as desired. So assume that color $M$ is present
on $G''$ (and not on $G'$).  Color $M-1$ must be present on $G'$, for
otherwise $\pch(G') \le M-2$.  If color $M-1$ is also used on $G''$, then
it again follows that ${\rm diam}(\widehat{G})>M-1$.  Hence we are left
with the situation that color $M$ is present on $G''$ and not on $G'$, while
$M-1$ is used on $G'$ and not on $G''$.  We now claim that the color $M-2$ is present in
both $G'$ and $G''$. For if this is not the case, then in any of $G'$ or $G''$
that is missing color $M-2$ relabeling all vertices colored with the highest color
by the color $M-2$ would yield an
$(M-2)$-packing coloring of $G'$ or $G''$, which is again not
possible. If $w',w''$ are the vertices in $G'$, resp.\ $G''$, with
$c(w')=c(w'')=M-2$, then $d_{\widehat{G}}(w',w'') \ge M-1$.  This in
turn implies  ${\rm diam}(\widehat{G})\ge M-1$, and so the claim
is proved.

Consider again vertices $w',w''$  in $G'$, resp.\ $G''$, with
$c(w')=c(w'')\ge M-2$. Since ${\rm diam}(G')\ge
d_{\widehat{G}}(w',x')$ and ${\rm diam}(G'')\ge
d_{\widehat{G}}(w'',x'')$, we infer that
\begin{eqnarray*}
2\,{\rm diam}(G')+1&=& {\rm diam}(G')+{\rm diam}(G'')+1\\
                 &\ge&
                 d_{\widehat{G}}(w',x')+d_{\widehat{G}}(w'',x'')+1
                 \\
                 &=& d_{\widehat{G}}(w',w'')\\
                 &\ge& M-1\,.
\end{eqnarray*}
Hence ${\rm diam}(G')\ge \lceil\frac{M}{2}\rceil-1$. Since clearly
${\rm diam}(G')\le {\rm diam}(G)+1$ holds, we conclude that
$${\rm diam}(G)\ge {\rm diam}(G')-1\ge \left\lceil\frac{M}{2}\right\rceil-2\,.$$
\qed

\begin{corollary}
\label{cor:14} There exists a cubic graph with packing chromatic
number larger than $13$.
\end{corollary}

\proof Let $G_{38}$ be the cubic graph of order $38$ with diameter
$4$ from \cite{afy-1986} shown in Figure~\ref{fig:G38}.

\begin{figure}[ht!]
\begin{center}
\begin{tikzpicture}[scale=4.0,style=thick]
    \vertex (0) at (0.97,.223) [scale=.75pt,fill=black]{};
    \vertex (1) at (0.9,.434) [scale=.75pt,fill=black]{};
    \vertex (2) at (0.782,.623) [scale=.75pt,fill=black]{};
    \vertex (3) at (0.623,.782) [scale=.75pt,fill=black]{};
    \vertex (4) at (0.434,.9) [scale=.75pt,fill=black]{};
    \vertex (5) at (0.223,.975) [scale=.75pt,fill=black]{};
    \vertex (6) at (0,1) [scale=.75pt,fill=black]{};
    \vertex (7) at (.97,-.023) [scale=.75pt,fill=black]{};
    \vertex (8) at (.9,-.234) [scale=.75pt,fill=black]{};
    \vertex (9) at (.782,-.423) [scale=.75pt,fill=black]{};
    \vertex (10) at (.623,-.582) [scale=.75pt,fill=black]{};
    \vertex (11) at (.434,-.7) [scale=.75pt,fill=black]{};
    \vertex (12) at (.18,-.775) [scale=.75pt,fill=black]{};
    \vertex (13) at (-.18,-.775) [scale=.75pt,fill=black]{};
    \vertex (14) at (-.434,-.7) [scale=.75pt,fill=black]{};
    \vertex (15) at (-.623,-.582) [scale=.75pt,fill=black]{};
    \vertex (16) at (-.782,-.423) [scale=.75pt,fill=black]{};
    \vertex (17) at (-.9,-.234) [scale=.75pt,fill=black]{};
    \vertex (18) at (-.97,-.023) [scale=.75pt,fill=black]{};
    \vertex (19) at (-.97,.223) [scale=.75pt,fill=black]{};
    \vertex (20) at (-.9,.434) [scale=.75pt,fill=black]{};
    \vertex (21) at (-.782,.623) [scale=.75pt,fill=black]{};
    \vertex (22) at (-.623,.782) [scale=.75pt,fill=black]{};
    \vertex (23) at (-.434,.9) [scale=.75pt,fill=black]{};
    \vertex (24) at (-.223,.975) [scale=.75pt,fill=black]{};
    \vertex (25) at (0,0.1) [scale=.75pt,fill=black]{};
    \vertex (26) at (0,0.3) [scale=.75pt,fill=black]{};
    \vertex (27) at (0,-.1) [scale=.75pt,fill=black]{};
    \vertex (28) at (-.143,.782) [scale=.75pt,fill=black]{};
    \vertex (29) at (.143,.782) [scale=.75pt,fill=black]{};
    \vertex (30) at (-.55,.55) [scale=.75pt,fill=black]{};
    \vertex (31) at (.55,.55) [scale=.75pt,fill=black]{};
    \vertex (32) at (0.75,.155) [scale=.75pt,fill=black]{};
    \vertex (33) at (-.75,.155) [scale=.75pt,fill=black]{};
    \vertex (34) at (-.57,-.3) [scale=.75pt,fill=black]{};
    \vertex (35) at (.57,-.3) [scale=.75pt,fill=black]{};
    \vertex (36) at (-.29,-.5) [scale=.75pt,fill=black]{};
    \vertex (37) at (.29,-.5) [scale=.75pt,fill=black]{};
    \path
        (0) edge[bend right=7] (1)
        (1) edge[bend right=7] (2)
        (2) edge[bend right=7] (3)
        (3) edge[bend right=7](4)
        (4) edge [bend right=7](5)
        (5) edge [bend right=7](6)
        (0) edge [bend left=7](7)
        (7) edge [bend left=7](8)
        (8) edge [bend left=7](9)
        (9) edge [bend left=7](10)
        (10) edge [bend left=7](11)
        (11) edge [bend left=7](12)
        (12) edge [bend left=7](13)
        (13) edge [bend left=7](14)
        (14) edge [bend left=7](15)
        (15) edge [bend left=7](16)
        (16) edge [bend left=7](17)
        (17) edge [bend left=7](18)
        (18) edge [bend left=7](19)
        (19) edge [bend left=7](20)
        (20) edge [bend left=7](21)
        (21) edge [bend left=7](22)
        (22) edge [bend left=7](23)
        (23) edge [bend left=7](24)
        (24) edge [bend left=7](6)
        (25) edge (26)
        (25) edge (27)
        (24) edge (28)
        (5) edge (29)
        (2) edge (31)
        (0) edge (32)
        (33) edge (19)
        (30) edge (21)
        (34) edge (16)
        (35) edge (9)
        (14) edge (36)
        (11) edge (37)
        (36) edge (32)
        (37) edge (33)
        (26) edge (20)
        (26) edge (1)
        (27) edge (15)
        (27) edge (10)
        (13) edge[bend right=11] (22)
        (12) edge[bend left=11] (3)
        (30) edge (31)
        (30) edge[bend right=10] (35)
        (31) edge[bend left=10] (34)
        (28) edge (34)
        (29) edge (35)
        (18) edge[bend right=5] (4)
        (7) edge[bend left=5] (23)
        (33) edge[bend right=15] (32)
        (37) edge[bend right=5] (28)
        (36) edge[bend left=5] (29)
        (17) edge[bend right=22] (8)
        (25) edge[bend left=16] (6)
        ;
\end{tikzpicture}
\caption{$G_{38}$}
\label{fig:G38}
\end{center}
\end{figure}

From~\cite[Proposition 6]{gt-2016} we know that
$\pch(G_ {38})=13$. We have checked by computer that
$\pch(S_e(G_{38}))=12$ holds for any edge $e$ of $G_{38}$. Assuming
that $M=13$ is the constant of Theorem~\ref{thm:diameter}, this
theorem implies that ${\rm diam}(G_{38})\ge
\lceil\frac{13}{2}\rceil-2=5$. However, since the diameter of
$G_{38}$ equals 4, we infer that $M$ cannot be 13. \qed

A closer look to the proof of Theorem~\ref{thm:diameter} reveals
that the graph constructed from two copies $G'_{38}$ and $G''_{38}$ of edge-subdivided
$G_{38}$ by connecting the vertices $x'$ and $x''$ is a graph of
order 78, say $G_{78}$ schematically shown in Figure~\ref{fig:G78}, such that $\pch(G_{78})\ge 14$.

\begin{figure}[ht!]
\begin{center}
\begin{tikzpicture}[scale=2.0,style=thick]
    \draw (0,0) ellipse (20pt and 12 pt);
    \draw (2,0) ellipse (20pt and 12 pt);
    \node (x1) at (0.32,0.02)[]{\small$x'$\normalsize};
    \draw (x1);
    \node (x2) at (1.72,0.02)[]{\small$x''$\normalsize};
    \draw (x1);
    \node (G1) at (-.2, 0.02)[]{$G_{38}'$};
    \draw(G1);
    \node (G2) at (2.2,0.02)[]{$G_{38}''$};
    \draw(G2);
    \vertex (0) at (0.45,0) [scale=.35pt,fill=black]{};
    \vertex (1) at (0.45,0.15) [scale=.35pt,fill=black]{};
    \vertex (2) at (0.45,-0.15) [scale=.35pt,fill=black]{};
    \vertex (3) at (1.55,-0.15) [scale=.35pt,fill=black]{};
    \vertex (4) at (1.55,0.15) [scale=.35pt,fill=black]{};
    \vertex (5) at (1.55,0) [scale=.35pt,fill=black]{};
    \path
        (0) edge (1)
        (0) edge (2)
        (0) edge (5)
        (5) edge (3)
        (5) edge (4);
\end{tikzpicture}
\caption{$G_{78}$}
\label{fig:G78}
\end{center}
\end{figure}

Motivated by the construction from the proof of
Theorem~\ref{thm:diameter}, we next consider what happens with the
packing chromatic number of an arbitrary graph when an edge is
subdivided.

\begin{theorem}
\label{thm:subdivision} For any graph $G$ with packing chromatic
number $j$,
\[\left\lfloor j/2 \right\rfloor +1 \le \chi_{\rho}(S_e(G)) \le
j+1.\] 
Moreover, for any $k\ge 2$ there exists a graph
$G$ with an edge $e$ such that $k=\pch(G)=\pch(S_e(G))-1$.
\end{theorem}
\proof Given a packing coloring $c$ of $G$, a packing coloring of
$S_e(G)$ can be obtained by using $c$ on vertices of $G$ and
coloring the new vertex with an additional color. Hence we get the upper
bound.

For the lower bound, let $G$ be a graph with $\pch(G)=j$ and consider any edge
$e=xy$ of $G$. Subdivide $e$ to get the graph $H=S_e(G)$.  That is,
we remove the edge $e$ from $G$ and replace it with the path
$x,z,y$.  Let $W_1, \ldots, W_r$ be an optimal packing coloring of
$H$.

We will construct a packing coloring of $G$.  Note that $\{x,y\}
\not\subseteq W_n$ for any $n\ge 2$. Fix $i$ such that $2 \le i \le
r$ and suppose there are vertices $u,v \in W_i$ such that $d_G(u,v)
= i$.  Since $W_i$ is an $i$-packing in $H$, we know that every
shortest $(u,v)$-path in $H$ contains the vertex $z$.  From among
all pairs of vertices in $W_i$ that are at distance $i$ in $G$
select $a_i,b_i$ such that $d_G(a_i,x)=t$ is the minimum and
$d_G(b_i,y)=s$, and so $t\le s$.  Thus
$d_H(a_i,b_i)=t+s+2=i+1$. Let $c\in W_i-\{a_i,b_i\}$. It now follows
that $d_G(c,x)> t$, for otherwise $d_H(c,a_i)\le 2t<i+1$, a
contradiction.  Similarly, $d_G(c,y)\ge s$, or otherwise it follows
that $d_H(a_i,c)=d_G(a_i,x)+2+d_G(y,c)<t+2+s=i+1$, again a
contradiction. For each such value of $i$ we remove the vertex $a_i$
from $W_i$ and place it into a set $X$ of vertices that will
eventually be ``recolored.'' For all pairs $u,v$ remaining in $W_i$
it follows that $d_G(u,v) \ge i+1$.  If $a_2$ as defined above
exists, then $a_2=x$.  It follows that $W_1$ is independent in $G$
and $|X|=m \le r-1$. Otherwise if $x$ and $y$ belong to $W_1$ place
vertex $x$ in the set $X$.  In this case $W_2$ is a $2$-packing in
$G$ and $|X|=m \le r-1$. Hence we can recolor the vertices in $X$
using colors $r+1,\ldots,r+m$ and this gives a packing coloring of
$G$ using at most $2r-1$ colors. That is, $\pch(G) \le 2r-1$.

 To prove the last assertion of the theorem, consider the
 following examples. For $k=2$, we have $2=\pch(P_3)=\pch(P_4)-1$, and
 $P_4=S_e(P_3)$. Let now $k\ge 3$. Recall~\cite[Lemma
 6]{bresar-2007} asserting that $\pch(S(K_k))=k+1$. Consider now the
 process of obtaining $S(K_k)$ from $K_k$ by subdividing each of the
 edges of $K_k$ one by one, and observe that in the beginning of this process
 $\pch(K_k)=k$, and at the end we have $\pch(S(K_k))=k+1$.
 Since in each step
 the packing chromatic number can increase by at most one, at some
 stage of the process we have graphs $G_i$ and $G_{i+1}$, such that
 $G_{i+1}=S_e(G_i)$ for some edge $e$ of $G_i$, and $\pch(G_i)=k$,
 $\pch(G_{i+1})=k+1$.
\qed

We do not know if the lower bound of
  Theorem~\ref{thm:subdivision} is sharp.
  On the other hand, it is possible that the subdivision of an edge
  decreases the packing chromatic number by 2. Consider the
  following examples. Let $n\ge 5$, and let $X_n$ be the graph
  obtained from the disjoint union of two copies of $K_n$, denoted by $U$ and $V$,
  by first joining a vertex $u$ of $U$ with a vertex $v$ of $V$,
  and then subdividing the edge $uv$ twice. Figure~\ref{fig:X5} depicts the graph $X_5$. Let $e=xy$ be the
  edge, where $x$ is adjacent to $u$, and $y$ is adjacent to $v$.

  \begin{figure}[h!]
\begin{center}
\begin{tikzpicture}[scale=2.0,style=thick]
    \vertex (0) at (0.3,0) [scale=.75pt,fill=black]{};
    \vertex (1) at (1.7,0) [scale=.75pt,fill=black]{};
    \vertex (2) at (.65,-.75) [scale=.75pt,fill=black]{};
    \vertex (3) at (1.35,-.75) [scale=.75pt,fill=black]{};
    \vertex (4) at (1,.55) [scale=.75pt,fill=black]{};
    \vertex (5) at (2.15,0) [label=below:$x$,scale=.75pt,fill=black]{};
    \vertex (6) at (2.6,0) [label=below:$y$,scale=.75pt,fill=black]{};
    \vertex (7) at (3.05,0) [scale=.75pt,fill=black]{};
    \vertex (8) at (4.45,0) [scale=.75pt,fill=black]{};
    \vertex (9) at (3.75,.55) [scale=.75pt,fill=black]{};
    \vertex (10) at (3.4,-.75) [scale=.75pt,fill=black]{};
    \vertex (11) at (4.1,-.75) [scale=.75pt,fill=black]{};

    \path
        (0) edge (1)
        (1) edge (2)
        (2) edge (3)
        (3) edge (4)
        (0) edge (2)
        (0) edge (3)
        (0) edge (4)
        (1) edge (3)
        (1) edge (4)
        (2) edge (4)
        (1) edge (5)
        (5) edge (6)
        (6) edge (7)
        (7) edge (8)
        (7) edge (9)
        (7) edge (10)
        (7) edge (11)
        (8) edge (9)
        (8) edge (10)
        (8) edge (11)
        (9) edge (10)
        (9) edge (11)
        (10) edge (11)
        ;

\end{tikzpicture}
\caption{$X_5$}
\label{fig:X5}
\end{center}
\end{figure}

  We claim that $\pch(X_n)=2n-3$. Let $c$ be an optimal packing
  coloring of $X_n$. Suppose first that $c(x)=1$ and $c(y)=2$.
Clearly, $c$ restricted to $U$ uses all the colors from $[n]$. On
the other hand, $c$ uses colors $1,3,4$ on $V$, while the other
$n-3$ vertices must receive new colors. Hence in this case, $c$ uses
$2n-3$ colors. Suppose next that one of the vertices $x$ and $y$ is
colored with a color $a$, where $a>2$. Then $c$ uses $n$ colors on
$U$, different from $a$, and $n-4$ new colors on $V$. Hence also in
this case $c$ uses $2n-3$ colors, which proves the claim.

Consider now the graph $S_e(X_n)$, and the following coloring $c$ of
this graph. Let $c(x)=c(y)=1$, and $c(v_{xy})=2$, where $v_{xy}$ is
the vertex obtained by subdividing the edge $xy$. The mapping $c$
restricted to $U$ uses colors from $[n]$ and restricted to $V$ uses
colors from $\{1,2,3,4,5\}$ together with $n-5$ new colors. Hence
$\pch(S_e(X_n))\le 2n-5$. The opposite inequality follows by the
observation that however $S_e(X_n)$ is colored, $n$ different colors
must be used on $U$, and out of these at most five colors can be
used also on $V$.  This shows that $\pch(S_e(X_n))= \pch(X_n)-2$.

Let us call an edge $e$ of graph $G$ {\em weak} if
$\pch(S_e(G))<\pch(G)-1$. We have not been able to find a graph $G$
such that all its edges are weak. We are inclined to believe that
there are no such graphs. From this point of view the following
consequence of Theorem~\ref{thm:diameter} is relevant.

\begin{corollary}
\label{cor:diameter} Suppose that there exists a constant $M$ such
that $\pch(H)\le M$ holds for any subcubic graph $H$, and let $G$ be
a subcubic graph such that $\pch(G)=M$. If there are no subcubic
graphs in which all edges are weak then $M\le 2\,{\rm diam}(G)+4$.
\end{corollary}

\section{Vertex deletion, edge deletion and contraction}
\label{sec:deletion}

Since the distances in a graph when an edge is removed can only
increase, it is clear that for any graph $G$, any vertex $v$ of $G$,
and any edge $e$ of $G$, we have
$$\pch(G-v)\le \pch(G)\quad {\rm and} \quad \pch(G-e)\le \pch(G).$$
On the other hand, there are no lower bounds for $\pch(G-v)$ and
$\pch(G-e)$. For the former operation, let $G_n$, $n\ge 4$, be the graph
obtained from the path $P_n$ by adding a vertex $x$ and making it
adjacent to all vertices of the path. Note that $\pch(G_n)\ge \lceil
\frac{n}{2}\rceil+1$, and since $G_n-x$ is isomorphic to $P_n$, we
have $\pch(G_n-x)=3$. To deal with  edge removal we state

\begin{proposition}
\label{prp:edgeremoval} For every positive integer $r$ there exists
a graph $G$ with an edge $e$ such that $\pch(G)-\pch(G-e)\ge r$.
\end{proposition}
\proof Consider the following construction. Let $k\ge 4$, and $n\ge
2k-2$. Let $A$ and $B$ be two copies of the graph $K_n$, and
$a,a'\in V(A)$, $b,b'\in V(B)$. The graph $G_{n,k}$ is obtained from
the disjoint union of $A$ and $B$ by connecting with an edge
vertices $a$ and $b$ and also connecting vertices $a'$ and $b'$, and
then replacing the edge $a'b'$ with a path of length $2k-1$. Figure~\ref{fig:G64} depicts the graph $G_{6,4}$.

  \begin{figure}[h!]
\begin{center}
\begin{tikzpicture}[scale=1.6,style=thick]
    \vertex (0) at (0.25,0) [scale=.75pt,fill=black]{};
    \vertex (1) at (0.25,-1) [scale=.75pt,fill=black]{};
    \vertex (2) at (1,-1.65) [scale=.75pt,fill=black]{};
    \vertex (3) at (1,.65) [scale=.75pt,fill=black]{};
    \vertex (4) at (1.75,0) [label=above:$a$,scale=.75pt,fill=black]{};
    \vertex (5) at (1.75,-1) [label=below:$a'$,scale=.75pt,fill=black]{};
    \vertex (6) at (2,-1) [scale=.75pt,fill=black]{};
    \vertex (7) at (2.25,-1) [scale=.75pt,fill=black]{};
    \vertex (8) at (2.5,-1) [scale=.75pt,fill=black]{};
    \vertex (9) at (2.75,-1) [scale=.75pt,fill=black]{};
    \vertex (10) at (3,-1) [scale=.75pt,fill=black]{};
    \vertex (11) at (3.25,-1) [scale=.75pt,fill=black]{};
    \vertex (12) at (3.5,-1) [label=below:$b'$,scale=.75pt,fill=black]{};
    \vertex (13) at (3.5,0) [label=above:$b$,scale=.75pt,fill=black]{};
    \vertex (14) at (4.25,.65) [scale=.75pt,fill=black]{};
    \vertex (15) at (4.25,-1.65) [scale=.75pt,fill=black]{};
    \vertex (16) at (5,0) [scale=.75pt,fill=black]{};
    \vertex (17) at (5,-1) [scale=.75pt,fill=black]{};

    \path
        (0) edge (1)
        (1) edge (2)
        (2) edge (3)
        (3) edge (4)
        (0) edge (2)
        (0) edge (3)
        (0) edge (4)
        (1) edge (3)
        (1) edge (4)
        (2) edge (4)
        (1) edge (5)
        (0) edge (5)
        (2) edge (5)
        (3) edge (5)
        (4) edge (5)
        (5) edge (6)
        (6) edge (7)
        (7) edge (8)
        (8) edge (9)
        (9) edge (10)
        (10) edge (11)
        (11) edge (12)
        (12) edge (13)
        (12) edge (14)
        (12) edge (15)
        (12) edge (16)
        (12) edge (17)
        (13) edge (14)
        (13) edge (15)
        (13) edge (16)
        (13) edge (17)
        (14) edge (15)
        (14) edge (16)
        (14) edge (17)
        (15) edge (16)
        (15) edge (17)
        (16) edge (17)
        (4) edge (16)

        ;

\end{tikzpicture}
\caption{$G_{6,4}$}
\label{fig:G64}
\end{center}
\end{figure}

We first claim that $\pch(G_{n,k})\ge 2n-2$. Note that $n$ colors
are used in any packing coloring on $A$. Since the distance between
a vertex of $A$ and a vertex of $B$ is at most $3$, we derive that
only colors $1$ and $2$ can be repeated in $B$, hence the claim.

Letting $G'_{n,k}=G_{n,k}-ab$ we next claim that $\pch(G'_{n,k})\le
2(n-k)+4$. Consider the following packing coloring of $G'_{n,k}$.
First color the path of length $2k-1$ between $a'$ and $b'$ with
colors from $\{1,2,3\}$. Because  in $G'_{n,k}$ every vertex in
$A\setminus\{a'\}$ is at distance $2k+1$ from any vertex in
$B\setminus\{b'\}$, we can use colors $4,\ldots,2k$ in both $A$ and
$B$. Note that this is possible because we have assumed that $n\ge
2k-2$, and hence the number of vertices in
$A\setminus\{a'\}$ and in $B\setminus\{{b}'\}$ is at least $2k-3$,
respectively. This in turn implies that the colors $4,\ldots,2k$ can
indeed be used twice. The remaining vertices are then colored by
unique colors. Consequently,
$$\pch(G'_{n,k}) \le 3+(2k-3)+[2(n-1)-2(2k-3)] = 2(n-k)+4\,.$$
It follows that $\pch(G_{n,k})-\pch(G'_{n,k})\ge
(2n-2)-[2(n-k)+4]=2k-6$. The assertion now follows. 
\qed


We next turn our attention to edge contractions. We denote the graph
obtained from $G$ by contracting its edge $e$ by $G|e$.

\begin{theorem}
\label{thm:contraction} If $G$ is a graph and $e$ an edge in $G$,
then $$\pch(G)-1\le \pch(G|e)\le 2\pch(G)\,.$$
\end{theorem}
\proof Let $e=xy$ be the edge that is contracted in a graph $G$, and
$v_{xy}$ the resulting vertex. For the lower bound, let $c$ be an
optimal packing coloring of $G|e$. We define the coloring $c'$ of
$G$ by letting $c'(x)=c(v_{xy})$, $c'(y)=\pch(G|e)+1$, and
$c'(z)=c(z)$ for any other vertex in $G$. Since the distances in $G$
are at least as large as the distances in $G|e$ between the
corresponding vertices, $c'$ is packing coloring of $G$. It follows
that $\pch(G)\le \pch(G|e)+1$.

For the upper bound let $c$ be an optimal packing coloring of $G$. We
define the coloring $c'$ of $G|e$ in two steps. First, let
$c'(v_{xy})=c(y)$, and $c'(z)=c(z)$ for any other vertex $z$ of
$G|e$. Let $i\in [\pch(G)]$, and let $x_i$ be a vertex of $G|e$ that
minimizes $d_{G|e}(z,v_{xy})$ over all $z\in V(G|e)$ with $c(z)=i$.
(Note that $x_i$ coincides with $v_{xy}$ for exactly one $i\in
[\pch(G)]$.) Then, in the second step, set $c'(x_i)=\pch(G)+i$. We
claim that $c'$ is a packing coloring of $G|e$.

Note that for any two vertices $a$ and $b$ of $G|e$ we have that
$d_{G|e}(a,b)$ is either $d_{G}(a,b)$ or $d_{G}(a,b)-1$. Moreover,
in the latter case there exists a shortest $(a,b)$-path in $G$ that
contains the edge $xy$. Suppose that there exist vertices $u$ and
$v$, both different from $x_i$, with $c'(u)=c'(v)=i$ such that
$d_{G|e}(u,v)=i$. Clearly, then in $G$ the edge $xy$ must lie on
some shortest $(u,v)$-path $P$ of length $i+1$. Hence we may assume
that $P$ is of the form $u-P'-x-y-P''-v$. We may also assume without
loss of generality that $d_G(x_i,x)\le d_G(x_i,y)$. Since $x_i$ is a
closest vertex to $v_{xy}$ among all vertices colored by $i$, we
derive that $d_G(x_i,x)<d_G(v,x)$, hence $d_G(u,x_i)\le
d_G(u,x)+d_G(x,x_i)< d_G(u,x)+d_G(x,v)= i+1$. This is a
contradiction with $c$ being a packing coloring of $G$, in which $u$
and $x_i$ are both colored by color $i$. This shows that
$c'$ is a packing coloring of $G|e$ with $2\pch(G)$ colors, hence
the proof of the upper bound is also complete.
\qed

Note that Theorem~\ref{thm:contraction} is in some sense dual to
Theorem~\ref{thm:subdivision}.
To see that the lower bound of Theorem~\ref{thm:contraction} is
sharp, just consider complete graphs. For the upper bound, similarly
as in Theorem~\ref{thm:subdivision}, we are not aware of any example
of a graph such that after the contraction of its edge the packing
chromatic number would increase by more than 2. On the other hand,
the graphs $S_e(X_n)$, as presented in Section~\ref{sec:subdivision}
show that the contraction of the edge $e$, yielding the graph $X_n$,
increases their packing chromatic number by $2$.

\section{Concluding remarks}
In this paper we answered a question of Gastineau and Togni~\cite{gt-2016} by
showing that there is a graph whose packing chromatic number is greater than $13$.
However, the problem from \cite{goddard-2008} concerning the existence of a constant
upper bound for the packing chromatic number on the class of cubic graphs remains
an interesting, unresolved problem.  It is possible that using
Theorem~\ref{thm:diameter} leads to subcubic graphs with increasing packing chromatic
number.  However, to prove this would require new methods since our approach in part
uses a computer.

Several open problems arise from considering local operations on
graphs and how these affect the packing chromatic number.  For
instance, the graph $G_{38}$ from Section~\ref{sec:subdivision}  has
the property that the subdivision of an arbitrary edge produces a
graph whose packing chromatic number is exactly one less than that
of $G_{38}$.  Cycles of the form $C_{4k+3}$ also share this
property. It would be interesting to know more about this class of
subdivision critical graphs. The examples $X_n$ from
Section~\ref{sec:subdivision} show that there exist graphs that have
an edge whose subdivision decreases the packing chromatic number by
2. As mentioned in Section~\ref{sec:subdivision} we
suspect that there does not exist a graph for which the subdivision of any of
its edges decreases the packing chromatic number by more then 1. In
other words, we suspect that there are no graphs with only weak
edges.

 Following the definition of graphs that are critical
with respect to ordinary chromatic number (i.e., the chromatic
number of any subgraph is less than that of the original graph) it
is natural to study graphs that are critical with respect to the
packing chromatic number. For graphs with no isolated vertices this
is equivalent to requiring that the packing chromatic number
decreases upon the removal of any edge.  Examples of these are
cycles whose order is not congruent to $0$ modulo $4$,  complete
graphs, and the Petersen graph.

\section*{Acknowledgements}
The authors wish to express their appreciation to Jernej Azarija for
the computations in the proof of Corollary~\ref{cor:14}. B.B. and
S.K. are supported in part by the Ministry of Science of Slovenia
under the grants P1-0297 and ARRS-BI-US/16-17-013. D.F.R. is supported by a grant from the
Simons Foundation (Grant Number 209654 to Douglas F. Rall).

\end{document}